\documentclass[12pt]{amsart}
\usepackage{amsmath}
\usepackage{amssymb}
\usepackage{amsfonts}
\usepackage{amsthm}
\usepackage{xcolor}
\usepackage[
    colorlinks=true,       
    linkcolor=blue,          
    citecolor=green,        
    filecolor=cyan,         
    urlcolor=magenta        
    ]{hyperref}
\usepackage{xurl}
\usepackage{mathrsfs}

\newtheorem{theo}{Theorem}[section]
\newtheorem*{theo*}{Theorem}

\newtheorem{lemm}[theo]{Lemma}

\newcommand{\N}{\mathbb{N}}

\newcommand{\C}{\mathbb{C}}
\newcommand{\lchi}{|L(1,\chi)|}

\title[]{A computable formula for evaluating the mean square sum of $L$-functions}
\begin{document}

 \author[N E Thomas]{Neha Elizabeth Thomas}
 \address{Department of Mathematics, University College, Thiruvananthapuram (Research Centre, University of Kerala), Kerala - 695034, India\\ORCID: 0000-0002-5473-5542}
 \email{nehathomas2009@gmail.com}
 \author[K V Namboothiri]{K Vishnu Namboothiri}
\address{Department of Mathematics, Baby John Memorial Government College, Chavara, Sankaramangalam, Kollam, Kerala - 691583, INDIA\\Department of Collegiate Education, Government of Kerala, India\\
ORCID: 0000-0003-4516-8117}
\email{kvnamboothiri@gmail.com}

 \begin{abstract}
  For Dirichlet characters $\chi$ mod $k$ where $k\geq 3$, we here give a computable formula for evaluating the mean square sums $\sum\limits_{\substack{\chi \text{ mod }k\\\chi(-1)=(-1)^r}}|L(r,\chi)|^2$ for any positive integer $r\geq 3$. We also give an inductive formula for computing the sum $\sum\limits_{\substack{1\leq m\leq k \\ (m, k)=1}}\frac{1}{\left(\sin\left(\frac{\pi m}{k}\right)\right)^{2n}}$ where $n$ is a positive integer in terms of Bernoulli numbers and binomial coefficients.
  \end{abstract}
 \keywords{$L$-functions, trigonometric sums, Jordan totient function, Euler totient function,
mean square sum, Chebyshev polynomials, Gauss sum, Ramanujan sum, Bernoulli numbers}
 \subjclass[2010]{11M06, 11L05, 11L03}

 \maketitle
\section{Introduction}
Let $k$ be a positive integer greater than or equal to 3. A Dirichlet character $\chi$ mod $k$ is defined to be \textit{odd} if $\chi(-1)=-1$ and \textit{even} if $\chi(-1)=1$. The \textit{Dirichlet $L$-function} $L(s,\chi)$ is defined by the infinite series $\sum\limits_{n=1}^{\infty}\frac{\chi(n)}{n^s}$, where $s\in\C$ with $\Re\,(s)>1$. It is an important function in number theory especially due to its connection with the Riemann zeta function $\zeta(s)$. For rational integer $r$, the problem of computing exact values of
\begin{align}\label{eqn:gen_sum}
 \sum\limits_{\substack{\chi \text{ mod }k\\\chi(-1)=(-1)^r}}|L(r,\chi)|^2
\end{align}
has been attempted in various cases by many.

 In 1982, Walum \cite{walum1982exact} gave an exact formula for the sum (\ref{eqn:gen_sum}) with $r=1$. Louboutin \cite{louboutin1999mean} computed  the sum of $\lchi^2$  over all odd primitive Dirichlet characters modulo $k$. See \cite[Chapter 6]{tom1976introduction} for the definition of primitivity of Dirichlet characters. In \cite{louboutin1993quelques}, Louboutin gave an exact formula for the sum of $\lchi^2$ over all odd Dirichlet characters in terms of the prime divisors of $k$ and the Euler totient function $\phi$.  In \cite{louboutin2001mean}  he derived an exact  formulae for the  general versions of these sums in two cases : $\chi$ even and $\chi$ odd. However, he did not suggest an efficient way to compute the coefficients appearing in his formula.

X. Lin \cite{lin2019mean} provided a general inductive formula for computing the sum (\ref{eqn:gen_sum}) which seems to be the only one computationally feasible formula existing for computing (\ref{eqn:gen_sum}) as of now. In \cite{okamoto2015various}, T. Okamoto and T. Onozuka  derived a formula for the sum
\begin{align*}
 \sum\limits_{\substack{\chi_1,\chi_2 \text{ mod }k\\\chi_1(-1)=(-1)^m\\\chi_2(-1)=(-1)^n}}L(m,\chi_1) L(n,\chi_2) L(m+n,\overline{\chi_1\chi_2})
\end{align*}
using the technique employed by S. Louboutin in \cite{louboutin2001mean}.  For more problems related to the sum of mean square values, please see   \cite{wu2012mean},  \cite{zhang2015hybrid}, \cite{okamoto2017mean}  and the more recent paper \cite{zhang2020certain}. Some other closely related problems were attempted by a few authors. See, for example, \cite{zhang2003problem}, \cite{ alkan2013averages} and \cite{ louboutin2015twisted}.

E. Alkan  derived an exact formula for the sums $\sum\limits_{\chi \, \text{odd}}|L(1,\chi)|^2$ and $\sum\limits_{\chi \, \text{even}}|L(2,\chi)|^2$  in terms of Jordan totient function and Euler totient function using weighted averages of Gauss and Ramanujan sums in \cite{alkan2011mean}. See the next section for the definition of these sums and functions. His starting point was the formula
 \begin{align*}                                                                                                                                                                   L(r,\chi)=\frac{i^r2^{r-1}\pi^r}{(-1)^{v+1}kr!}\sum\limits_{q=0}^{2[\frac{r}{2}]}\binom{r}{q}
 B_q S(r-q,\chi) \text{ \cite[Theorem 1]{alkan2011values}},
 \end{align*}
 where $v$ and $r \geq 1$ have same parity, $B_q$ are the Bernoulli numbers and $S(r-q, \chi)=\sum\limits_{j=1}^k \left(\frac{j}{k}\right)^{r-q}G(z,\chi)$. Here $G(z,\chi)$ denotes the Gauss sum.  The most difficult and complicated part of his evaluation was computing sums of the form
 $$  \sum\limits_{\substack{1\leq m\leq k \\ (m, k)=1}}\left(\sum\limits_{j=1}^{k-1}j^pe^{\frac{2\pi imj}{k}}\right)\left(\sum\limits_{s=1}^{k-1}s^qe^{\frac{2\pi ims}{k}}\right)$$ and $$\sum\limits_{\substack{1\leq m\leq k \\ (m, k)=1}}\left(\sum\limits_{j=1}^{k-1}j^pe^{\frac{2\pi imj}{k}}\right)\left(\sum\limits_{s=1}^{k-1}s^qe^{-\frac{2\pi ims}{k}}\right)$$ for which he did not suggest any general technique, but simplified it only for his special case.

Note that Alkan's techniques in \cite{alkan2011mean} were completely different from the other derivations which were employed to compute the square sums. It was also observed by Alkan that his techniques could be used to determine the sum (\ref{eqn:gen_sum}) for any larger $r$, but at the cost of increasingly complex computations. We followed his techniques and gave exact formulae for $\sum\limits_{\chi(-1)=(-1)^r}|L(r,\chi)|^2$ when $r=3$ and $r=4$ in \cite{thomas2023mean}.

Using the similar techniques as that of Alkan in \cite{alkan2011mean}  and these authors in \cite{thomas2023mean}, we here give a computationally feasible general formula for evaluating (\ref{eqn:gen_sum}) in two cases: $r$ odd (Theorem \ref{identity:oddcasefinalformula}), and $r$ even (Theorem
\ref{identity:evencasefinalformula}). In both these cases, the formulae are almost the same, and involve the computation of  certain binomial coefficients and sums of the form \begin{align*}\sum\limits_{\substack{1\leq m\leq k \\ (m, k)=1}}\frac{1}{\left(\sin\left(\frac{\pi m}{k}\right)\right)^{2n}},\,n = 1,2,\ldots,r.
\end{align*}
 For finding these sine sums, we give an inductive formula (Lemma \ref{lem:sineevenpowersum}) which only involves the evaluation of certain binomial coefficients.

 \section{Notations and basic results}
In our computations to determine $\sum\limits_{\substack{\chi(\text{mod }  k)\\\chi(-1)=(-1)^r}}|L(r, \chi)|^2$, we use certain properties of Chebyshev polynomials of first and second kinds. A detailed account of these polynomials and their properties can be found in \cite[Section 13.3]{arfken2005mathematical}. The \textit{Chebyshev polynomials of the first kind} $T_n$ has the representation
\begin{align*}
T_n(x)=\frac{n}{2}\sum\limits_{k=0}^{\lfloor\frac{n}{2}\rfloor}(-1)^k\frac{(n-k-1)!}{k!(n-2k)!}(2x)^{n-2k}\text{ \cite[Equation 13.119a]{arfken2005mathematical}}
\end{align*}
and the \textit{second kind} $U_n$ has the representation
\begin{align*}
U_n(x)=\sum\limits_{k=0}^{\lfloor\frac{n}{2}\rfloor}(-1)^k\frac{(n-k)!}{k!(n-2k)!}(2x)^{n-2k}\text{ \cite[Equation 13.119b]{arfken2005mathematical}},
\end{align*}
where $\lfloor x \rfloor$ denote the greatest integer not exceeding $x$. Now these polynomials can be used to relate the sine and cosine of $\theta$ and $n\theta$ via the identities
$T_n(\cos\theta)=\cos (n\theta)$ and $U_n(\cos\theta)\sin\theta=\sin((n+1)\theta)$ \cite[Equations 13.114a and 13.115a]{arfken2005mathematical}. 

The \textit{Gauss sum} $G(z,\chi)$ for any complex number $z$ is defined by
 \begin{align*}
 G(z,\chi):= \sum\limits_{m=1}^{k}\chi(m)e^{\frac{2\pi imz}{k}}.
 \end{align*}

 The \textit{Jordan totient function} $J_k(n)$ is defined by  \begin{align*}
J_k(n) := n^k\prod_{\substack{p|n\\p\text{ prime}}}\left(1-\frac{1}{p^k}\right).
\end{align*}

For $n\geq 0$,  the $n$th Bernoulli number is denoted by $B_n$. It is known that $B_0=1$ and $B_1 =-\frac{1}{2}$. For any other odd positive integer $n$, we have $B_n=0$. For the other properties of Bernoulli numbers, please see \cite[Chapter 12]{tom1976introduction}.

From \cite{alkan2011values} we have the following identities.
 For positive integer $n$, we have \begin{align}\label{alk1}
 \sum\limits_{j=1}^{k-1}j^ne^{\frac{2\pi imj}{k}}=\sum\limits_{j=1}^{n}\binom{n}{j}k^j\lim \limits_{w\rightarrow\frac{2\pi im}{k}}\frac{d^{n-j}}{dw^{n-j}}\left(\frac{1}{e^w-1}\right).
 \end{align}
 
 For any $q\geq1$, we have 
 \begin{align}\label{alk2}
 \frac{d^{q}}{dw^{q}}\left(\frac{1}{e^w-1}\right)=\sum\limits_{j=1}^{q+1}\frac{A_{q,j}}{(e^w-1)^j},
 \end{align} where
 \begin{align}\label{aqj}
 A_{q,j}:=\sum\limits_{r=0}^{j-1}(-1)^{r+q}\binom{j-1}{r}(j-r)^q
 \end{align}
 for $1\leq j\leq q+1$.
 Combining the two Equations (\ref{alk1}) and (\ref{alk2}) we get
 \begin{align}\label{alk12}
 \sum\limits_{j=1}^{k-1}j^ne^{\frac{2\pi imj}{k}}=\sum\limits_{j=1}^{n}\binom{n}{j}k^j\sum\limits_{\alpha=1}^{n-j+1}\frac{A_{n-j,\alpha}}{(e^{\frac{2\pi im}{k}}-1)^\alpha}.
 \end{align}
Since the sum $\sum\limits_{\substack{1\leq m\leq k \\ (m, k)=1}}\frac{1}{\left(\sin\left(\frac{\pi m}{k}\right)\right)^{n}}$ occurs at many places in this paper, for notational brevity, we will use the notation $SIN(n)$ to replace it.
\section{Main Results and proofs}\label{main-results-section}
We start with the evaluation of the real part of 
 \begin{align*}
  \sum\limits_{\substack{1\leq m\leq k \\ (m, k)=1}}\left(\sum\limits_{j=1}^{k-1}j^pe^{\frac{2\pi imj}{k}}\right)\left(\sum\limits_{s=1}^{k-1}s^qe^{\frac{2\pi ims}{k}}\right).
  \end{align*}
We use Equation (\ref{alk12}) to evaluate the above and get it to be equal to
 \begin{align}\label{sumj^ps^q}
  &\sum\limits_{\substack{1\leq m\leq k \\ (m, k)=1}}\sum\limits_{j=1}^{k-1}j^pe^{\frac{2\pi imj}{k}}\sum\limits_{s=1}^{k-1}s^qe^{\frac{2\pi ims}{k}}\\
  &= \sum\limits_{\substack{1\leq m\leq k \\ (m, k)=1}}\sum\limits_{s=1}^{k-1}s^qe^{\frac{2\pi ims}{k}}\sum\limits_{j=1}^{p}\binom{p}{j}k^j\sum\limits_{\alpha=1}^{p-j+1}\frac{A_{p-j, \alpha}}{\left(e^{\frac{2\pi im}{k}}-1\right)^\alpha}\nonumber\\
  &=\sum\limits_{\substack{1\leq m\leq k \\ (m, k)=1}}\sum\limits_{j=1}^{p}\binom{p}{j}k^j\sum\limits_{\alpha=1}^{p-j+1}\frac{A_{p-j, \alpha}}{\left(e^{\frac{2\pi im}{k}}-1\right)^\alpha}\sum\limits_{s=1}^{k-1}s^qe^{\frac{2\pi ims}{k}}\nonumber\\
  &=\sum\limits_{\substack{1\leq m\leq k \\ (m, k)=1}}\sum\limits_{j=1}^{p}\binom{p}{j}k^j\sum\limits_{\alpha=1}^{p-j+1}\frac{A_{p-j, \alpha}}{\left(e^{\frac{2\pi im}{k}}-1\right)^\alpha}\left(\sum\limits_{s=1}^{q}\binom{q}{s}k^s\sum\limits_{\beta=1}^{q-s+1}\frac{A_{q-s, \beta}}{\left(e^{\frac{2\pi im}{k}}-1\right)^\beta}\right)\nonumber\\
  &=\sum\limits_{j=1}^{p}\binom{p}{j}\sum\limits_{s=1}^{q}\binom{q}{s}k^{j+s}\sum\limits_{\substack{1\leq m\leq k \\ (m, k)=1}}\sum\limits_{\alpha=1}^{p-j+1}\sum\limits_{\beta=1}^{q-s+1}\frac{A_{p-j, \alpha}A_{q-s, \beta}}{\left(e^{\frac{2\pi im}{k}}-1\right)^{\alpha+\beta}}.\nonumber
\end{align}
Now by the Euler's formula $e^{i\theta} = \cos\theta + i \sin \theta$, we have
\begin{align}\label{1/e^w-1}
\frac{1}{\left(e^{\frac{2\pi im}{k}}-1\right)^{\alpha+\beta}}&=\frac{\left(e^{-\frac{\pi i m}{k}}\right)^{\alpha+\beta}}{\left(e^{\frac{\pi i m}{k}}-e^{-\frac{\pi i m}{k}}\right)^{\alpha+\beta}}\\
&=\frac{\cos\left(\left(\frac{\pi m}{k}\right)(\alpha+\beta)\right)-i\sin\left(\left(\frac{\pi m}{k}\right)(\alpha+\beta)\right)}{(2i)^{\alpha+\beta}\left(\sin\left(\frac{\pi m}{k}\right)\right)^{\alpha+\beta}}.\nonumber
\end{align}
Since we are interested only in the real part of the left-hand side of Equation (\ref{sumj^ps^q}), only  the real part of $\frac{1}{\left(e^{\frac{2\pi im}{k}}-1\right)^{\alpha+\beta}}$ in Equation (\ref{1/e^w-1}) is relevant to us. To compute it, we consider two cases determined by the parity of $\alpha+\beta$ as follows.\\
Case 1: $\alpha+\beta$   even

 When $\alpha+\beta$ is even, $i^{\alpha+\beta}$ is real and only the term  $\frac{\cos\left(\left(\frac{\pi m}{k}\right)(\alpha+\beta)\right)}{(2i)^{\alpha+\beta}\left(\sin\left(\frac{\pi m}{k}\right)\right)^{\alpha+\beta}}$ appears in the real part of Equation (\ref{1/e^w-1}). Now
 \begin{align}\label{Chebyshev}
\cos\left(\left(\frac{\pi m}{k}\right)(\alpha+\beta)\right)&=T_{\alpha+\beta}\left(cos\left(\frac{\pi m}{k}\right)\right) \\
&=\frac{\alpha+\beta}{2}\sum\limits_{c=0}^{\lfloor\frac{\alpha+\beta}{2}\rfloor}(-1)^c\frac{(\alpha+\beta-c-1)!}{c!(\alpha+\beta-2c)!}\left(2\cos\left(\frac{\pi m}{k}\right)\right)^{\alpha+\beta-2c}.\nonumber
 \end{align}
 Hence
 \begin{align*}
 &\frac{\cos\left(\left(\frac{\pi m}{k}\right)(\alpha+\beta)\right)}{(2i)^{\alpha+\beta}\left(\sin\left(\frac{\pi m}{k}\right)\right)^{\alpha+\beta}}\\
&=\frac{\frac{\alpha+\beta}{2}\sum\limits_{c=0}^{\lfloor\frac{\alpha+\beta}{2}\rfloor}(-1)^c\frac{(\alpha+\beta-c-1)!}{c!(\alpha+\beta-2c)!}\left(2\cos\left(\frac{\pi m}{k}\right)\right)^{\alpha+\beta-2c}}{2^{\alpha+\beta}i^{\alpha+\beta}\left(\sin\left(\frac{\pi m}{k}\right)\right)^{\alpha+\beta}}\\
&=(\alpha+\beta)\sum\limits_{c=0}^{\lfloor\frac{\alpha+\beta}{2}\rfloor}\frac{(-1)^{c+\frac{\alpha+\beta}{2}}(\alpha+\beta-c-1)!}{2^{2c+1}c!(\alpha+\beta-2c)!\left(\sin\left(\frac{\pi m}{k}\right)\right)^{\alpha+\beta}}\left(\cos\left(\frac{\pi m}{k}\right)\right)^{\alpha+\beta-2c}.
 \end{align*}
 Now
 \begin{align*}
 \left(\cos\left(\frac{\pi m}{k}\right)\right)^{\alpha+\beta-2c}&= \left(1-\sin^2\left(\frac{\pi m}{k}\right)\right)^{\frac{\alpha+\beta-2c}{2}}\\&=\sum\limits_{d=0}^{\frac{\alpha+\beta-2c}{2}}(-1)^d\binom{\frac{\alpha+\beta-2c}{2}}{d}\left(\sin\left(\frac{\pi m}{k}\right)\right)^{2d}.
 \end{align*}
Therefore
\begin{align}\label{alpha+betaeven}
&\frac{\cos\left(\left(\frac{\pi m}{k}\right)(\alpha+\beta)\right)}{(2i)^{\alpha+\beta}\left(\sin\left(\frac{\pi m}{k}\right)\right)^{\alpha+\beta}}\\
=&(\alpha+\beta)\sum\limits_{c=0}^{\frac{\alpha+\beta}{2}}\frac{(-1)^{c+\frac{\alpha+\beta}{2}}(\alpha+\beta-c-1)!}{2^{2c+1}c!(\alpha+\beta-2c)!}
\sum\limits_{d=0}^{\frac{\alpha+\beta-2c}{2}}\frac{(-1)^d\binom{\frac{\alpha+\beta-2c}{2}}{d}}{\left(\sin\left(\frac{\pi m}{k}\right)\right)^{\alpha+\beta-2d}}.\nonumber
\end{align} 
Case 2: $\alpha+\beta$   odd

 When $\alpha+\beta$ is odd, $i^{\alpha+\beta}$ is imaginary and the term  $-\frac{\sin\left(\left(\frac{\pi m}{k}\right)(\alpha+\beta)\right)}{2^{\alpha+\beta}i^{\alpha+\beta-1}\left(\sin\left(\frac{\pi m}{k}\right)\right)^{\alpha+\beta}}$ appears in the real part of Equation (\ref{1/e^w-1}). Now
 \begin{align*}
&\sin\left(\left(\frac{\pi m}{k}\right)(\alpha+\beta)\right)\\&=U_{\alpha+\beta-1}\left(\cos\left(\frac{\pi m}{k}\right)\right) \sin\left(\frac{\pi m}{k}\right)\\
&=\left[\sum\limits_{e=0}^{\lfloor\frac{\alpha+\beta-1}{2}\rfloor}(-1)^e\frac{(\alpha+\beta-e-1)!}{e!(\alpha+\beta-2e-1)!}\left(2\cos\left(\frac{\pi m}{k}\right)\right)^{\alpha+\beta-2e-1}\right]\sin\left(\frac{\pi m}{k}\right).
 \end{align*}
Also
\begin{align*}
\left(\cos\left(\frac{\pi m}{k}\right)\right)^{\alpha+\beta-2e-1}&= \left(1-\sin^2\left(\frac{\pi m}{k}\right)\right)^{\frac{\alpha+\beta-2e-1}{2}}\\&=\sum\limits_{f=0}^{\frac{\alpha+\beta-2c-1}{2}}(-1)^f\binom{\frac{\alpha+\beta-2e-1}{2}}{f}\left(\sin\left(\frac{\pi m}{k}\right)\right)^{2f}.
\end{align*}
Hence
\begin{align}\label{abetaodd}
&-\frac{\sin\left(\left(\frac{\pi m}{k}\right)(\alpha+\beta)\right)}{2^{\alpha+\beta}i^{\alpha+\beta-1}\left(\sin\left(\frac{\pi m}{k}\right)\right)^{\alpha+\beta}}=(-1)^{\frac{\alpha+\beta+1}{2}}\frac{\sin\left(\left(\frac{\pi m}{k}\right)(\alpha+\beta)\right)}{2^{\alpha+\beta}\left(\sin\left(\frac{\pi m}{k}\right)\right)^{\alpha+\beta}}\\
&=\sum\limits_{e=0}^{\frac{\alpha+\beta-1}{2}}\frac{(-1)^{e+\frac{\alpha+\beta+1}{2}}(\alpha+\beta-e-1)!}{2^{2e+1}e!(\alpha+\beta-2e-1)!}\sum\limits_{f=0}^{\frac{\alpha+\beta-2e-1}{2}}\frac{(-1)^f\binom{\frac{\alpha+\beta-2e-1}{2}}{f}}{\left(\sin\left(\frac{\pi m}{k}\right)\right)^{\alpha+\beta-2f-1}}.\nonumber
\end{align}

With these, from Equation (\ref{sumj^ps^q}), we have

\begin{lemm}\label{realjs}
\begin{align*}
&\Re \left(\sum\limits_{\substack{1\leq m\leq k \\ (m, k)=1}}\left(\sum\limits_{j=1}^{k-1}j^pe^{\frac{2\pi imj}{k}}\right)\left(\sum\limits_{s=1}^{k-1}s^qe^{\frac{2\pi ims}{k}}\right)\right)\\ \nonumber
=&\sum\limits_{j=1}^{p}\sum\limits_{s=1}^{q}\binom{p}{j}\binom{q}{s}k^{j+s}\sum\limits_{\alpha=1}^{p-j+1}A_{p-j, \alpha}\Bigg[\Bigg(\sum\limits_{\substack{\beta=1\\ \alpha+\beta \text{ even}}}^{q-s+1}A_{q-s, \beta}\sum\limits_{c=0}^{\frac{\alpha+\beta}{2}}(\alpha+\beta)\\\nonumber
&\times\frac{(-1)^{c+\frac{\alpha+\beta}{2}}(\alpha+\beta-c-1)!}{2^{2c+1}c!(\alpha+\beta-2c)!}
\sum\limits_{d=0}^{\frac{\alpha+\beta-2c}{2}}(-1)^d\binom{\frac{\alpha+\beta-2c}{2}}{d}SIN(\alpha+\beta-2d)\Bigg)\\\nonumber
&+\Bigg(\sum\limits_{\substack{\beta=1\\ \alpha+\beta \text{ odd}}}^{q-s+1}A_{q-s, \beta}\sum\limits_{e=0}^{\frac{\alpha+\beta-1}{2}}\frac{(-1)^{e+\frac{\alpha+\beta+1}{2}}(\alpha+\beta-e-1)!}{2^{2e+1}e!(\alpha+\beta-2e-1)!}\\\nonumber
&\times\sum\limits_{f=0}^{\frac{\alpha+\beta-2e-1}{2}}(-1)^f\binom{\frac{\alpha+\beta-2e-1}{2}}{f}SIN(\alpha+\beta-2f-1)\Bigg)\Bigg].\nonumber
\end{align*}
\end{lemm}

Note that from the computations we performed to arrive at the result above, a similar looking identity can be computed for \begin{align*}
\Re \left(\sum\limits_{\substack{1\leq m\leq k \\ (m, k)=1}}\left(\sum\limits_{j=1}^{k-1}j^pe^{\frac{2\pi imj}{k}}\right)\right)                                                                                                               \end{align*}

 \subsection{Computing the sum when $r$ is odd}
 Now we compute $\sum\limits_{\chi(\text{mod }  k)}|L(r, \chi)|^2$ over odd characters $\chi$ when $r$ is odd. We will use the  identity in Lemma \ref{realjs} in these computations.
 By \cite[Theorem 1]{alkan2011values}, we have   
 \begin{align}\label{1}
 L(r,\chi)=\frac{i^r2^{r-1}\pi^r}{(-1)^{v+1}kr!}\sum\limits_{q=0}^{2[\frac{r}{2}]}\binom{r}{q}
 B_q S(r-q,\chi),
 \end{align}
 where $ v$  and  $r\geq1$ have same parity, and so in this case, $(-1)^{v+1}=1$.

 Here
\begin{align*}
S(r-q, \chi)=\sum\limits_{j=1}^k\left(\frac{j}{k}\right)^{r-q}G(j, \chi)=\sum\limits_{j=1}^k \left(\frac{j}{k}\right)^{r-q}\sum\limits_{m=1}^{k-1}\chi(m)e^{\frac{2\pi imj}{k}}.
\end{align*}
 Since $r$ is odd, write $r=2h+1$, where $h=0,1,\cdots$.
Thus we have 
\begin{align*}
L(2h+1, \chi)=C_1 \sum \limits_{q=0}^{2h}\binom{2h+1}{q}B_q \sum\limits_{j=1}^{k}\left(\frac{j}{k}\right)^{2h+1-q}\sum\limits_{m=1}^{k-1}\chi(m)e^{\frac{2\pi imj}{k}},
\end{align*}
where $C_1=\frac{i^{2h+1}2^{2h}\pi^{2h+1}}{k(2h+1)!}$. Hence
\begin{align*}
\sum\limits_{\substack{\chi(\text{mod }  k) \\ \chi \text{ odd}}}&|L(2h+1, \chi)|^2\\=&\sum\limits_{\substack{\chi(\text{mod }  k) \\ \chi \text{ odd}}}L(2h+1, \chi)\overline{L(2h+1, \chi)}\nonumber\\
=&C_1\overline{C_1}\sum\limits_{q_1=0}^{2h}\sum\limits_{q_2=0}^{2h}B_{q_1}B_{q_2}\binom{2h+1}{q_1}\binom{2h+1}{q_2}\sum\limits_{j=1}^{k-1}\left(\frac{j}{k}\right)^{2h+1-q_1}\sum\limits_{s=1}^{k-1}\left(\frac{s}{k}\right)^{2h+1-q_2}\nonumber\\ &\times\sum\limits_{m=1}^{k-1}\sum\limits_{n=1}^{k-1}e^{\frac{2\pi i(mj-ns)}{k}}\sum\limits_{\substack{\chi(\text{mod }  k) \\ \chi \text{ odd}}}\chi(m)\overline{\chi(n)}.\nonumber
\end{align*}
There are exactly $\frac{\phi(k)}{2}$ odd Dirichlet characters modulo $k\geq3$. Therefore $\sum\limits_{\substack{\chi(\text{mod }  k) \\ \chi \text{ odd}}}\chi(u)=\frac{\phi(k)}{2}$ or $-\frac{\phi(k)}{2}$ depending on if $u=1$ or $u=-1$ respectively. If $u$ is neither of these modulo $k$, then  $\sum\limits_{\substack{\chi(\text{mod }  k) \\ \chi \text{ odd}}}\chi(u)=0$. Thus
 $\sum\limits_{\substack{\chi(\text{mod }  k) \\ \chi \text{ odd}}}\chi(m)\overline{\chi(n)}=\sum\limits_{\substack{\chi(\text{mod }  k) \\ \chi \text{ odd}}}\chi(mn^{-1})=\frac{\phi(k)}{2} \text{ or } -\frac{\phi(k)}{2}$\\
when $m=n$ or $m=k-n$ respectively and  
 $\sum\limits_{\substack{\chi(\text{mod }  k) \\ \chi \text{ odd}}}\chi(m)\overline{\chi(n)}=0 $
 otherwise.  So we get
 \begin{align}\label{prod}
 \sum\limits_{\substack{\chi(\text{mod }  k) \\ \chi \text{ odd}}}&|L(2h+1, \chi)|^2\nonumber\\
= &|C_1|^2\frac{\phi(k)}{2}\sum\limits_{q_1=0}^{2h}\sum\limits_{q_2=0}^{2h}B_{q_1}B_{q_2}\frac{1}{k^{4h+2-q_1-q_2}}\binom{2h+1}{q_1}\binom{2h+1}{q_2}\nonumber \\
&\times\Bigg[\sum\limits_{j=1}^{k-1} \sum\limits_{s=1}^{k-1}j^{2h+1-q_1}s^{2h+1-q_2}\sum\limits_{\substack{1\leq m\leq k \\ (m, k)=1}}e^{\frac{2\pi im(j-s)}{k}}\nonumber\\
&-\sum\limits_{j=1}^{k-1} \sum\limits_{s=1}^{k-1}j^{2h+1-q_1}s^{2h+1-q_2}\sum\limits_{\substack{1\leq m\leq k \\ (m, k)=1}}e^{\frac{2\pi im(j+s)}{k}}\Bigg]\nonumber\\
=&|C_1|^2\frac{\phi(k)}{2}\sum\limits_{q_1=0}^{2h}\sum\limits_{q_2=0}^{2h}B_{q_1}B_{q_2}\frac{1}{k^{4h+2-q_1-q_2}}\binom{2h+1}{q_1}\binom{2h+1}{q_2}\\
&\times\Bigg[\sum\limits_{\substack{1\leq m\leq k \\ (m, k)=1}}\left(\sum\limits_{j=1}^{k-1}j^{2h+1-q_1}e^{\frac{2\pi imj}{k}}\right)\left(\sum\limits_{s=1}^{k-1}s^{2h+1-q_2}e^{-\frac{2\pi ims}{k}}\right)\nonumber\\
&-\sum\limits_{\substack{1\leq m\leq k \\ (m, k)=1}}\left(\sum\limits_{j=1}^{k-1}j^{2h+1-q_1}e^{\frac{2\pi imj}{k}}\right)\left(\sum\limits_{s=1}^{k-1}s^{2h+1-q_2}e^{\frac{2\pi ims}{k}}\right)\Bigg].\nonumber
\end{align}

We use identity \cite[Proposition 2.1]{thomas2023exponential} to proceed further. It states that if $f(q)=\sum\limits_{s=1}^{k-1}s^qe^{\frac{-2\pi ims}{k}}$ and $g(p)= \sum\limits_{s=1}^{k-1}s^{p}e^{\frac{2\pi ims}{k}}$, then
$f(q) =-k^q + \sum\limits_{a=0}^{q-1}(-1)^{q-a}\binom{q}{a}k^{a}g(q-a).$
So 
\begin{align*}
\sum\limits_{\substack{1\leq m\leq k \\ (m, k)=1}}g(p)f(q)&=\sum\limits_{\substack{1\leq m\leq k \\ (m, k)=1}}g(p)\left(-k^q + \sum\limits_{a=0}^{q-1}(-1)^{q-a}\binom{q}{a}k^{a}g(q-a)\right)\\
&=-k^q\sum\limits_{\substack{1\leq m\leq k \\ (m, k)=1}}g(p) + \sum\limits_{a=0}^{q-1}(-1)^{q-a}\binom{q}{a}k^{a}\sum\limits_{\substack{1\leq m\leq k \\ (m, k)=1}}g(p)g(q-a).
\end{align*}
Thus
\begin{align}\label{g(p)f(q)}
&\sum\limits_{\substack{1\leq m\leq k \\ (m, k)=1}}\left(\sum\limits_{j=1}^{k-1}j^pe^{\frac{2\pi imj}{k}}\right)\left(\sum\limits_{s=1}^{k-1}s^qe^{-\frac{2\pi ims}{k}}\right)\\
=&-k^q\sum\limits_{\substack{1\leq m\leq k \\ (m, k)=1}}\left(\sum\limits_{j=1}^{k-1}j^pe^{\frac{2\pi imj}{k}}\right)+\nonumber\\ &\sum\limits_{a=0}^{q-1}(-1)^{q-a}\binom{q}{a}k^{a}\sum\limits_{\substack{1\leq m\leq k \\ (m, k)=1}}\left(\sum\limits_{j=1}^{k-1}j^pe^{\frac{2\pi imj}{k}}\right)\left(\sum\limits_{s=1}^{k-1}s^{q-a}e^{\frac{2\pi ims}{k}}\right).\nonumber
\end{align}
We use Equation \eqref{g(p)f(q)} with $p=2h+1-q_1$ and $q = 2h+1-q_2$ in Equation (\ref{prod})  to get
\begin{align}\label{producttermfinal}
&\sum\limits_{\substack{\chi(\text{mod }  k) \\ \chi \text{ odd}}}|L(2h+1, \chi)|^2 =|C_1|^2\frac{\phi(k)}{2}\\
&\Bigg[\sum\limits_{q_1=0}^{2h}\sum\limits_{q_2=0}^{2h}B_{q_1}B_{q_2}\frac{1}{k^{4h+2-q_1-q_2}}\binom{2h+1}{q_1}\binom{2h+1}{q_2}\times -k^{2h+1-q_2}\nonumber\\
&\sum\limits_{\substack{1\leq m\leq k \\ (m, k)=1}}\left(\sum\limits_{j=1}^{k-1}j^{2h+1-q_1}e^{\frac{2\pi imj}{k}}\right)+\nonumber\\
&\sum\limits_{q_1=0}^{2h}\sum\limits_{q_2=0}^{2h}B_{q_1}B_{q_2}\binom{2h+1}{q_1}\binom{2h+1}{q_2}\sum\limits_{a=0}^{2h-q_2}(-1)^{2h+1-q_2-a}
\times \binom{2h+1-q_2}{a}\nonumber\\
&
k^{a-4h-2+q_1+q_2}\sum\limits_{\substack{1\leq m\leq k \\ (m, k)=1}}\left(\sum\limits_{j=1}^{k-1}j^{2h+1-q_1}e^{\frac{2\pi imj}{k}}\right)\left(\sum\limits_{s=1}^{k-1}s^{2h+1-q_2-a}e^{\frac{2\pi ims}{k}}\right)\nonumber\\
-&\sum\limits_{q_1=0}^{2h}\sum\limits_{q_2=0}^{2h}B_{q_1}B_{q_2}k^{-4h-2+q_1+q_2}\binom{2h+1}{q_1}\binom{2h+1}{q_2}\nonumber\\
&\sum\limits_{\substack{1\leq m\leq k \\ (m, k)=1}}\left(\sum\limits_{j=1}^{k-1}j^{2h+1-q_1}e^{\frac{2\pi imj}{k}}\right)\left(\sum\limits_{s=1}^{k-1}s^{2h+1-q_2}e^{\frac{2\pi ims}{k}}\right)\Bigg].\nonumber
\end{align}

In the first sum

\begin{align*}
 \Sigma_0:=&\sum\limits_{q_1=0}^{2h}\sum\limits_{q_2=0}^{2h}B_{q_1}B_{q_2}\frac{1}{k^{4h+2-q_1-q_2}}\binom{2h+1}{q_1}\binom{2h+1}{q_2}\times -k^{2h+1-q_2}\\
&\sum\limits_{\substack{1\leq m\leq k \\ (m, k)=1}}\left(\sum\limits_{j=1}^{k-1}j^{2h+1-q_1}e^{\frac{2\pi imj}{k}}\right)
\end{align*}
in the right-hand side of the above, if we rearrange the summations in such a way that the innermost sum becomes $\sum\limits_{q_2=0}^{2h}B_{q_2}\binom{2h+1}{q_2}$, then by \cite[Theorem 12.15]{tom1976introduction}, if $2h+1\geq 3$, we have  $0=B_{2h+1}=\sum\limits_{q_2=0}^{2h+1}B_{q_2}\binom{2h+1}{q_2}=\sum\limits_{q_2=0}^{2h}B_{q_2}\binom{2h+1}{q_2}$. So $\Sigma_0$ is equal to 0 in this case. When $2h+1=1$, from  \cite[Equations 2.13-2.15]{alkan2011mean} we get
\begin{align*}
 \sum\limits_{\substack{1\leq m\leq k \\ (m, k)=1}}\left(\sum\limits_{j=1}^{k-1}j^{2h+1-q_1}e^{\frac{2\pi imj}{k}}\right)=-\frac{k\phi(k)}{2}.
\end{align*}
So $\Sigma_0=\frac{\phi(k)}{2}$ when $2h+1=1$.

Also, in Equation (\ref{producttermfinal}), since the left-hand side is a real number, we need to look only at the real part of the right-hand side. The right-hand side consists of two summations (the second and third sums) over the variables $q_1$ and $q_2$. The imaginary part in both these summations come from the sumations over $1\leq m\leq k, (m,k)=1$ in which the exponential function is present. For further computations, we consider two summations over $q_1$ and $q_2$ seperately which we denote as

\begin{align}
 \Sigma_1 := &\sum\limits_{q_1=0}^{2h}\sum\limits_{q_2=0}^{2h}B_{q_1}B_{q_2}\binom{2h+1}{q_1}\binom{2h+1}{q_2}\sum\limits_{a=0}^{2h-q_2}(-1)^{2h+1-q_2-a}
\times \binom{2h+1-q_2}{a}\label{sigma5_original}\\
&k^{a-4h-2+q_1+q_2}\Re \left(\sum\limits_{\substack{1\leq m\leq k \\ (m, k)=1}}\left(\sum\limits_{j=1}^{k-1}j^{2h+1-q_1}e^{\frac{2\pi imj}{k}}\right)\left(\sum\limits_{s=1}^{k-1}s^{2h+1-q_2-a}e^{\frac{2\pi ims}{k}}\right)\right)\nonumber
\end{align}
and
\begin{align*}
 \Sigma_2 := \sum\limits_{q_1=0}^{2h}\sum\limits_{q_2=0}^{2h}&B_{q_1}B_{q_2}\binom{2h+1}{q_1}\binom{2h+1}{q_2}k^{-4h-2+q_1+q_2}\\
&\Re \left(\sum\limits_{\substack{1\leq m\leq k \\ (m, k)=1}}\left(\sum\limits_{j=1}^{k-1}j^{2h+1-q_1}e^{\frac{2\pi imj}{k}}\right)\left(\sum\limits_{s=1}^{k-1}s^{2h+1-q_2}e^{\frac{2\pi ims}{k}}\right)\right)\nonumber
\end{align*}
so that when $2h+1\geq 3$, Equation (\ref{producttermfinal}) can be rewritten as
\begin{align*}
 \sum\limits_{\substack{\chi(\text{mod }  k) \\ \chi \text{ odd}}}|L(2h+1, \chi)|^2 =|C_1|^2\frac{\phi(k)}{2}[\Sigma_1-\Sigma_2].
 \end{align*}

Note that by Lemma \ref{realjs}  we already have expanded expressions for the real parts included in the above two sums $\Sigma_1$ and $\Sigma_2$. We use them to get

\begin{align}\label{sigma1}
  \Sigma_1 =& \sum\limits_{q_1=0}^{2h}\sum\limits_{q_2=0}^{2h}\sum\limits_{a=0}^{2h-q_2}\sum\limits_{j=1}^{2h+1-q_1}\sum\limits_{s=1}^{2h+1-q_2-a}\sum\limits_{\alpha=1}^{2h+2-q_1-j}B_{q_1}B_{q_2}\binom{2h+1}{q_1}\binom{2h+1}{q_2}
 \binom{2h+1-q_2}{a}\\\nonumber
 &\times \binom{2h+1-q_1}{j}\binom{2h+1-q_2-a}{s}k^{a-4h-2+q_1+q_2+j+s} A_{2h+1-q_1-j, \alpha}\\\nonumber
 &\times\Bigg[\Bigg(\sum\limits_{\substack{\beta=1\\ \alpha+\beta \text{ even}}}^{2h+2-q_2-a-s} \sum\limits_{c=0}^{\frac{\alpha+\beta}{2}}\sum\limits_{d=0}^{\frac{\alpha+\beta-2c}{2}}(-1)^{2h+1-q_2-a+c+d+\frac{\alpha+\beta}{2}}A_{2h+1-q_2-a-s, \beta}\\\nonumber
&\times\frac{(\alpha+\beta)(\alpha+\beta-c-1)!}{2^{2c+1}c!(\alpha+\beta-2c)!}
\binom{\frac{\alpha+\beta-2c}{2}}{d}SIN(\alpha+\beta-2d)\Bigg)\\\nonumber
& +\Bigg(\sum\limits_{\substack{\beta=1\\ \alpha+\beta \text{ odd}}}^{2h+2-q_2-a-s}\sum\limits_{e=0}^{\frac{\alpha+\beta-1}{2}}\sum\limits_{f=0}^{\frac{\alpha+\beta-2e-1}{2}}(-1)^{2h+1-q_2-a+e+f+\frac{\alpha+\beta+1}{2}}A_{2h+1-q_2-a-s, \beta}\\\nonumber
&\times\frac{(\alpha+\beta-e-1)!}{2^{2e+1}e!(\alpha+\beta-2e-1)!}\binom{\frac{\alpha+\beta-2e-1}{2}}{f}SIN(\alpha+\beta-2f-1)\Bigg)\Bigg]\nonumber
\end{align}
and
\begin{align}\label{sigma2}
  \Sigma_2 =& \sum\limits_{q_1=0}^{2h}\sum\limits_{q_2=0}^{2h}\sum\limits_{j=1}^{2h+1-q_1}\sum\limits_{s=1}^{2h+1-q_2}\sum\limits_{\alpha=1}^{2h+2-q_1-j}B_{q_1}B_{q_2}\binom{2h+1}{q_1}\binom{2h+1}{q_2}
\\&\times
\binom{2h+1-q_1}{j}\binom{2h+1-q_2}{s}k^{-4h-2+q_1+q_2+j+s} A_{2h+1-q_1-j, \alpha} \nonumber\\&\times\Bigg[\Bigg(\sum\limits_{\substack{\beta=1\\ \alpha+\beta \text{ even}}}^{2h+2-q_2-s} \sum\limits_{c=0}^{\frac{\alpha+\beta}{2}}\sum\limits_{d=0}^{\frac{\alpha+\beta-2c}{2}}(-1)^{c+d+\frac{\alpha+\beta}{2}}A_{2h+1-q_2-s, \beta}\nonumber\\
&\times\frac{(\alpha+\beta)(\alpha+\beta-c-1)!}{2^{2c+1}c!(\alpha+\beta-2c)!}
\binom{\frac{\alpha+\beta-2c}{2}}{d}SIN(\alpha+\beta-2d)\Bigg)\nonumber\\
 &+\Bigg(\sum\limits_{\substack{\beta=1\\ \alpha+\beta \text{ odd}}}^{2h+2-q_2-s}\sum\limits_{e=0}^{\frac{\alpha+\beta-1}{2}}\sum\limits_{f=0}^{\frac{\alpha+\beta-2e-1}{2}}(-1)^{e+f+\frac{\alpha+\beta+1}{2}}A_{2h+1-q_2-s, \beta}\nonumber\\
&\times\frac{(\alpha+\beta-e-1)!}{2^{2e+1}e!(\alpha+\beta-2e-1)!}\binom{\frac{\alpha+\beta-2e-1}{2}}{f}SIN(\alpha+\beta-2f-1)\Bigg)\Bigg],\nonumber
\end{align}
where $A_{q,j}$ is given by Equation (\ref{aqj}). Now we show that $\Sigma_1=-\Sigma_2$ or equivalently $\Sigma_1+\Sigma_2=0$. Suppose we fix the values of the variables $q_1=q_0, j=j_0, s=s_0, \alpha=\alpha_0,\beta=\beta_0$ (with $\alpha+\beta$ even), $c=c_0,d=d_0$, then we get a term in $\Sigma_1+\Sigma_2$ as

\begin{align*}
 &\sum\limits_{q_2=0}^{2h}\sum\limits_{a=0}^{2h-q_2}B_{q_0}B_{q_2}\binom{2h+1}{q_0}\binom{2h+1}{q_2}
 \binom{2h+1-q_2}{a}\\&\times
\binom{2h+1-q_0}{j_0}\binom{2h+1-q_2-a}{s_0}k^{a-4h-2+q_0+q_2+j_0+s_0} A_{2h+1-q_0-j_0, \alpha_0}\\&\times\Bigg[(-1)^{2h+1-q_2-a+c_0+d_0+\frac{\alpha_0+\beta_0}{2}}A_{2h+1-q_2-a-s_0, \beta_0}\\
&\times\frac{(\alpha_0+\beta_0)(\alpha_0+\beta_0-c_0-1)!}{2^{2c_0+1}c_0!(\alpha_0+\beta_0-2c_0)!}
\binom{\frac{\alpha_0+\beta_0-2c_0}{2}}{d_0}SIN(\alpha_0+\beta_0-2d_0)\Bigg]
\end{align*}
\begin{align*}
&+\sum\limits_{q_2=0}^{2h}B_{q_0}B_{q_2}\binom{2h+1}{q_0}\binom{2h+1}{q_2}
 \\&\times
\binom{2h+1-q_0}{j_0}\binom{2h+1-q_2}{s_0}k^{-4h-2+q_0+q_2+j_0+s_0} A_{2h+1-q_0-j_0, \alpha_0}\\&\times\Bigg[(-1)^{c_0+d_0+\frac{\alpha_0+\beta_0}{2}}A_{2h+1-q_2-s_0, \beta_0}\\
&\times\frac{(\alpha_0+\beta_0)(\alpha_0+\beta_0-c_0-1)!}{2^{2c_0+1}c_0!(\alpha_0+\beta_0-2c_0)!}\binom{\frac{\alpha_0+\beta_0-2c_0}{2}}{d_0}SIN(\alpha_0+\beta_0-2d_0)\Bigg].
\end{align*}

In the second sum above over $q_2$, fix $q_2=m$ where $0\leq m\leq 2h$. Also, in the first sum, fix $a+q_2=m$. Then the above sum becomes

\begin{align}\label{expression}
 &\sum\limits_{q_2=0}^{2h}B_{q_0}B_{q_2}\binom{2h+1}{q_0}\binom{2h+1}{q_2}
 \binom{2h+1-q_2}{m-q_2}\\&\times
\binom{2h+1-q_0}{j_0}\binom{2h+1-m}{s_0}k^{m-4h-2+q_0+j_0+s_0} A_{2h+1-q_0-j_0, \alpha_0}\nonumber\\&\times\Bigg[(-1)^{2h+1-m+c_0+d_0+\frac{\alpha_0+\beta_0}{2}}A_{2h+1-m-s_0, \beta_0}\nonumber\\
&\times\frac{(\alpha_0+\beta_0)(\alpha_0+\beta_0-c_0-1)!}{2^{2c_0+1}c_0!(\alpha_0+\beta_0-2c_0)!}
\binom{\frac{\alpha_0+\beta_0-2c_0}{2}}{d_0}SIN(\alpha_0+\beta_0-2d_0)\Bigg]\nonumber
\end{align}
\begin{align*}
&+\nonumber\\
&B_{q_0}B_{m}\binom{2h+1}{q_0}\binom{2h+1}{m}\nonumber
 \\&\times
\binom{2h+1-q_0}{j_0}\binom{2h+1-m}{s_0}k^{-4h-2+q_0+m+j_0+s_0} A_{2h+1-q_0-j_0, \alpha_0}\nonumber\\&\times\Bigg[(-1)^{c_0+d_0+\frac{\alpha_0+\beta_0}{2}}A_{2h+1-m-s_0, \beta_0}\nonumber\\
&\times
\frac{(\alpha_0+\beta_0)(\alpha_0+\beta_0-c_0-1)!}{2^{2c_0+1}c_0!(\alpha_0+\beta_0-2c_0)!}\binom{\frac{\alpha_0+\beta_0-2c_0}{2}}{d_0}SIN(\alpha_0+\beta_0-2d_0)\Bigg]\nonumber.
\end{align*}

Now we sum the expression (\ref{expression}) over $m$ where $0\leq m\leq 2h$. Note that, since $a+q_2=m$ and $a$ and $q_2$ are both nonnegative, $q_2\leq m$. Also, $(-1)^{2h+1-m}=(-1)^{m+1}$. On summing over $m$, after collecting similar terms and making some rearrangements of the sums, we get

\begin{align}\label{expression1}
 &B_{q_0}\binom{2h+1}{q_0}\binom{2h+1-q_0}{j_0}\binom{2h+1-m}{s_0}k^{m-4h-2+q_0+j_0+s_0} \\
 &A_{2h+1-q_0-j_0, \alpha_0} A_{2h+1-m-s_0, \beta_0} (-1)^{c_0+d_0+\frac{\alpha_0+\beta_0}{2}}\nonumber\\&\times\frac{(\alpha_0+\beta_0)(\alpha_0+\beta_0-c_0-1)!}{2^{2c_0+1}c_0!(\alpha_0+\beta_0-2c_0)!}
\binom{\frac{\alpha_0+\beta_0-2c_0}{2}}{d_0}SIN(\alpha_0+\beta_0-2d_0)\nonumber
\\
 &\times\sum\limits_{m=0}^{2h}\left((-1)^{m+1}\sum\limits_{q_2=0}^{m}B_{q_2}\binom{2h+1}{q_2}
 \binom{2h+1-q_2}{m-q_2} +B_{m}\binom{2h+1}{m}\right)\nonumber.
\end{align}

Now $\binom{2h+1}{q_2} \binom{2h+1-q_2}{m-q_2} = \binom{2h+1}{m}  \binom{m}{q_2}$. Also, by \cite[Theorem 12.15]{tom1976introduction}, we have $\sum\limits_{q_2=0}^{m}B_{q_2}\binom{m}{q_2}=B_m$ for $m\geq 2$.
So
\begin{align*}
 &\sum\limits_{m=0}^{2h}\left((-1)^{m+1}\sum\limits_{q_2=0}^{m}B_{q_2}\binom{2h+1}{q_2}
 \binom{2h+1-q_2}{m-q_2} +B_{m}\binom{2h+1}{m}\right)\\
 &=\sum\limits_{m=0}^{2h}\binom{2h+1}{m}\left((-1)^{m+1}\sum\limits_{q_2=0}^{m}B_{q_2}\binom{m}{q_2} + B_{m}\right)\nonumber.
\end{align*}

Now \begin{align*}
     (-1)^{m+1}\sum\limits_{q_2=0}^{m}B_{q_2}\binom{m}{q_2} + B_{m}=\begin{cases}
                                                                     B_m+B_m=0\text{ when }m\geq2 \text{ and odd}\\
                                                                     -B_m+B_m=0\text{ when }m\geq2 \text{ and even}.
                                                                    \end{cases}    \end{align*}

When $m=0$ or 1, with the values $B_0=1$ and $B_1=-\frac{1}{2}$, then also we get
\begin{align*}
     (-1)^{m+1}\sum\limits_{q_2=0}^{m}B_{q_2}\binom{m}{q_2} + B_{m}=0.
     \end{align*}
     Therefore expression (\ref{expression1}) and hence expression (\ref{expression}) evaluates to $0$.
     We can proceed in the similar way to arrive at the same conclusion when $\alpha+\beta$ is odd.
     Hence it follows that $\Sigma_1+\Sigma_2 = 0$  or $\Sigma_1=-\Sigma_2$.

  Thus

 \begin{theo}\label{identity:oddcasefinalformula}
 For $h \geq 1, h\in \N$ we have
\begin{align*}
 \sum\limits_{\substack{\chi(\text{mod }  k) \\ \chi \text{ odd}}}|L(2h+1, \chi)|^2 =&-\frac{2^{4h}\pi^{4h+2}}{((2h+1)!)^2k^2}\phi(k)\Sigma_2,
\end{align*}
where $\Sigma_2$ is as in Equation (\ref{sigma2}).
\end{theo}
In fact, we can evaluate $\sum\limits_{\substack{\chi(\text{mod }  k) \\ \chi \text{ odd}}}|L(1, \chi)|^2$ also using the above formula; we just need to add the term $|C_1|^2\frac{\phi(k)}{2}\Sigma_0=\frac{\pi^2(\phi(k))^2}{4k^2}$ to the right-hand side of the above for that.

In the sum $\Sigma_2$, all the sumations are easy to evaluate using some computational software like \texttt{SageMath} except the sum $SIN(r)=\sum\limits_{\substack{1\leq m\leq k \\ (m, k)=1}}\frac{1}{\left(\sin\left(\frac{\pi m}{k}\right)\right)^{r}}$. Note also that exponent of $\sin\left(\frac{\pi m}{k}\right)$ is even in $\Sigma_2$. Now we proceed to give an inductive formula for computing $SIN(r)$ when $r$ is even. By \cite[Theorems 1 and 2]{alkan2011values} we have
\begin{align*}
L(r,\chi)=C_2\sum\limits_{m=1}^{k-1}\chi(m)\sum\limits_{q=0}^{2[\frac{r}{2}]}\binom{r}{q}\frac{B_q}{k^{r-q}}\sum\limits_{j=1}^{r-q}\binom{r-q}{j}k^j\sum\limits_{\alpha=1}^{r-q-j+1}\frac{A_{r-q-j, \alpha}}{\left(e^{\frac{2\pi im}{k}}-1\right)^\alpha},
\end{align*}
where $C_2=\frac{i^r2^{r-1}\pi^r}{(-1)^{v+1}kr!}$ with $r$ and $v$ are of same parity so that $(-1)^{v+1}=-1$ and $i^r=(-1)^{\frac{r}{2}}$ when $r$ is even. When $\chi=\chi_0$, we have $2[\frac{r}{2}]=r$ and so
\begin{align}\label{lrchinot}
 L(r,\chi_0)=C_2\sum\limits_{q=0}^{r}\binom{r}{q}\frac{B_q}{k^{r-q}}\sum\limits_{j=1}^{r-q}\binom{r-q}{j}k^j\sum\limits_{\alpha=1}^{r-q-j+1}\sum\limits_{\substack{1\leq m\leq k \\ (m, k)=1}}\frac{A_{r-q-j, \alpha}}{\left(e^{\frac{2\pi im}{k}}-1\right)^\alpha}.
\end{align}
But we also have
$L(r, \chi_0)=\zeta(r)\prod\limits_{\substack{p|k \\ p \text{ prime}}}\left(1-\frac{1}{p^{r}}\right)=\zeta(r)\frac{J_{r}(k)}{k^{r}}$ \cite[Theorem 11.7]{tom1976introduction} and $\zeta(r)=\frac{(-1)^{\frac{r}{2}+1}(2\pi)^{r}B_{r}}{2r!}$ \cite[Theorem 12.17]{tom1976introduction}. On combining both these, we get
\begin{align}\label{lrchinotapostol}
 L(r, \chi_0)=\frac{(-1)^{\frac{r}{2}+1}(2\pi)^{r}B_{r}}{2r!}\frac{J_{r}(k)}{k^{r}}.
 \end{align}
 Hence $L(r, \chi_0)$ is real when $r$ is even. In the right-hand side of the Equation (\ref{lrchinot}), the imaginary part is contributed by the term $\frac{1}{\left(e^{\frac{2\pi im}{k}}-1\right)^{\alpha}}$. We need to consider only the real part of this term to evaluate  $L(r, \chi_0)$ in Equation (\ref{lrchinot}). As in the beginning of section \ref{main-results-section}, corresponding to the two cases $\alpha$  even and  $\alpha$ odd, the real parts of $\frac{1}{\left(e^{\frac{2\pi im}{k}}-1\right)^{\alpha}}$ are different. When $\alpha$ is even, from Equation (\ref{alpha+betaeven}), we get

\begin{align*}
\Re\left(\frac{1}{\left(e^{\frac{2\pi im}{k}}-1\right)^{\alpha}}\right)&= \frac{\cos\left(\left(\frac{\pi m}{k}\right)\alpha\right)}{(2i)^{\alpha}\left(\sin\left(\frac{\pi m}{k}\right)\right)^{\alpha}} \\
&=\alpha\sum\limits_{c=0}^{\frac{\alpha}{2}}\frac{(-1)^{c+\frac{\alpha}{2}}(\alpha-c-1)!}{2^{2c+1}c!(\alpha-2c)!}
\sum\limits_{d=0}^{\frac{\alpha-2c}{2}}\frac{(-1)^d\binom{\frac{\alpha-2c}{2}}{d}}{\left(\sin\left(\frac{\pi m}{k}\right)\right)^{\alpha-2d}}
\end{align*} and when $\alpha$ is odd, from Equation (\ref{abetaodd}), we get

\begin{align*}
\Re\left(\frac{1}{\left(e^{\frac{2\pi im}{k}}-1\right)^{\alpha}}\right)&=-\frac{\sin\left(\left(\frac{\pi m}{k}\right)\alpha\right)}{(2i)^{\alpha-1}\left(\sin\left(\frac{\pi m}{k}\right)\right)^{\alpha}}\\
&=\sum\limits_{e=0}^{\frac{\alpha-1}{2}}\frac{(-1)^{e+\frac{\alpha+1}{2}}(\alpha-1-e)!}{2^{2e+1}e!(\alpha-1-2e)!}\sum\limits_{f=0}^{\frac{\alpha-1-2c}{2}}\frac{(-1)^f\binom{\frac{\alpha-1-2e}{2}}{f}}{\left(\sin\left(\frac{\pi m}{k}\right)\right)^{\alpha-2f-1}}.
\end{align*}
Thus Equation (\ref{lrchinot}) becomes
\begin{align}\label{L(r,x_0)}
 &L(r,\chi_0)  \\=& C_2\sum\limits_{q=0}^{r}\binom{r}{q}\frac{B_q}{k^{r-q}}\sum\limits_{j=1}^{r-q}\binom{r-q}{j}k^j\Bigg[\Bigg(\sum\limits_{\substack{\alpha=1\\ \alpha \text{ even}}}^{r-q-j+1}A_{r-q-j, \alpha}\nonumber\\&\times\alpha\sum\limits_{c=0}^{\frac{\alpha}{2}}\frac{(-1)^{c+\frac{\alpha}{2}}(\alpha-c-1)!}{2^{2c+1}c!(\alpha-2c)!}
\sum\limits_{d=0}^{\frac{\alpha-2c}{2}}(-1)^d\binom{\frac{\alpha-2c}{2}}{d}SIN(\alpha-2d)\Bigg)\nonumber
\end{align}
\begin{align*}
&+\Bigg(\sum\limits_{\substack{\alpha=1\\ \alpha \text{ odd}}}^{r-q-j+1}A_{r-q-j, \alpha}\sum\limits_{e=0}^{\frac{\alpha-1}{2}}\frac{(-1)^{e+\frac{\alpha+1}{2}}(\alpha-1-e)!}{2^{2e+1}e!(\alpha-1-2e)!}\nonumber\\&\times\sum\limits_{f=0}^{\frac{\alpha-1-2e}{2}}(-1)^f\binom{\frac{\alpha-1-2e}{2}}{f}SIN(\alpha-2f-1)\Bigg)\Bigg].\nonumber
\end{align*}

 In $SIN(\alpha-2d) = \sum\limits_{\substack{1\leq m\leq k \\ (m, k)=1}}\frac{1}{\left(\sin\left(\frac{\pi m}{k}\right)\right)^{\alpha-2d}}$ and $SIN(\alpha-2f-1) = \sum\limits_{\substack{1\leq m\leq k \\ (m, k)=1}}\frac{1}{\left(\sin\left(\frac{\pi m}{k}\right)\right)^{\alpha-2f-1}}$ in the right-hand side of the above, the largest exponent of $\frac{1}{\sin\left(\frac{\pi m}{k}\right)}$is $r$ which occurs when $\alpha-2d=r$. All the other exponents of $\frac{1}{\sin\left(\frac{\pi m}{k}\right)}$ are less than $r$. Note that $\alpha-2d=r$ when $\alpha$ is maximum and $d$ is minimum. That is, when $q+j$ is minimum (which is 1) and $d=0$. We move the term consisting of $SIN(r)$ from  the right-hand side of Equation (\ref{L(r,x_0)}) to the left-hand side and substitute the value of $L(r,\chi_0)$ obtained in Equation (\ref{lrchinotapostol}) to get
\begin{align*}
 &
 \Bigg(\frac{C_2r}{k^{r-1}}A_{r-1, r}r\sum\limits_{c=0}^{\frac{r}{2}}\frac{(-1)^{c+\frac{r}{2}}(r-c-1)!}{2^{2c+1}c!(r-2c)!}\Bigg)
SIN(r)\nonumber\\&=\frac{(-1)^{\frac{r}{2}+1}(2\pi)^{r}B_{r}}{2r!}\frac{J_{r}(k)}{k^{r}}-C_2\sum\limits_{q=0}^{r}\binom{r}{q}\frac{B_q}{k^{r-q}}\sum\limits_{j=1}^{r-q}\binom{r-q}{j}k^j\Bigg[\Bigg(\sum\limits_{\substack{\alpha=1\\ \alpha \text{ even}}}^{r-q-j+1}A_{r-q-j, \alpha}\nonumber\\&\times\alpha\sum\limits_{c=0}^{\frac{\alpha}{2}}\frac{(-1)^{c+\frac{\alpha}{2}}(\alpha-c-1)!}{2^{2c+1}c!(\alpha-2c)!}
\sum\limits_{\substack{d=0\\ \alpha-2d\neq r}}^{\frac{\alpha-2c}{2}}(-1)^d\binom{\frac{\alpha-2c}{2}}{d}SIN(\alpha-2d)\Bigg)\nonumber\\
&+\Bigg(\sum\limits_{\substack{\alpha=1\\ \alpha \text{ odd}}}^{r-q-j+1}A_{r-q-j, \alpha}\sum\limits_{e=0}^{\frac{\alpha-1}{2}}\frac{(-1)^{\frac{\alpha+1}{2}+e}(\alpha-1-e)!}{2^{2e+1}e!(\alpha-1-2e)!}\nonumber\\&\times\sum\limits_{f=0}^{\frac{\alpha-1-2e}{2}}(-1)^f\binom{\frac{\alpha-1-2e}{2}}{f}SIN(\alpha-2f-1)\Bigg)\Bigg].\nonumber
\end{align*}

We also have $A_{r-1, r} = (-1)^{r-1}(r-1)!$ \cite[Theorem 2]{alkan2011values}. When $\alpha +\beta = r$ and $\frac{m}{k}=2$ in Equation (\ref{Chebyshev}), we get $r\sum\limits_{c=0}^{\frac{r}{2}}\frac{(-1)^{c+\frac{r}{2}}(r-c-1)!}{2^{2c+1}c!(r-2c)!}= \frac{(-1)^{\frac{r}{2}}}{2^r}$. Substituting this and the value of $C_2$, we get the coefficient of $SIN(r)$ in the above equation as $\frac{\pi^r}{2k^r}$. On further simplification of constants we obtain
\begin{lemm}\label{lem:sineevenpowersum}
For an even positive integer $r$, we have
\begin{align*}
 &SIN(r)\\&= \frac{(-1)^{\frac{r}{2}+1}2^{r}B_{r}}{r!}J_{r}(k)+\frac{(-1)^{\frac{r}{2}}2^{r}k^{r-1}}{r!}\sum\limits_{q=0}^{r}\binom{r}{q}\frac{B_q}{k^{r-q}}\sum\limits_{j=1}^{r-q}\binom{r-q}{j}k^j\Bigg[\sum\limits_{\alpha=1}^{r-q-j+1}A_{r-q-j, \alpha}\nonumber\\&\times E(\alpha)\sum\limits_{c=0}^{[\frac{\alpha}{2}]}\frac{(-1)^{c+[\frac{\alpha+1}{2}]}(\alpha-c-1)!}{2^{2c+1}c!(2[\frac{\alpha}{2}]-2c)!}
\sum\limits_{\substack{d=0\\ \alpha-2d\neq r}}^{[\frac{\alpha}{2}]-c}(-1)^d\binom{[\frac{\alpha}{2}]-c}{d}SIN\left(2\left[\frac{\alpha}{2}\right]-2d\right)\Bigg],\nonumber
\end{align*}
where $E(\alpha)=
\begin{cases}
\alpha \text{ if } \alpha \text{ is even}\\
1\text{ if } \alpha \text{ is odd}
  \end{cases}$.
\end{lemm}
The above result together with the fact that $SIN(0)=\phi(k)$ can be used to evaluate $SIN(r)$ for any even positive integer $r$.  This completes the evaluation of $\sum\limits_{\substack{\chi(\text{mod }  k) \\ \chi \text{ odd}}}|L(2h+1, \chi)|^2 $ using the expression given in Theorem \ref{identity:oddcasefinalformula}.

\subsection{Computing the sum when $r$ is even}

Now we proceed to give a similar formula for  $\sum\limits_{\substack{\chi(\text{mod }  k) \\ \chi \text{ even}}}|L(r, \chi)|^2$ when $r\geq 4$ is even. Let $r=2h$ where $h\geq2$. Note that
\begin{align*}
S(0,\chi)=\sum\limits_{j=1}^{k}\sum\limits_{m=1}^{k-1}\chi(m)e^{\frac{2\pi imj}{k}}=\sum\limits_{m=1}^{k-1}\chi(m)\sum\limits_{j=1}^{k}e^{\frac{2\pi imj}{k}}=0.
 \end{align*}
 Also for $h\geq2$, we have $B_{2h-1}=0$.
Hence Equation (\ref{1}) with $r=2h$ gives
\begin{align*}
  L(2h,\chi)&=C_3 \sum \limits_{q=0}^{2h-2}\binom{2h}{q}B_q S(2h-q,\chi)\\
  &=C_3 \sum \limits_{q=0}^{2h-2}\binom{2h}{q}B_q \sum\limits_{j=1}^{k-1}\left(\frac{j}{k}\right)^{2h-q}\sum\limits_{m=1}^{k-1}\chi(m)e^{\frac{2\pi imj}{k}},\nonumber
\end{align*}
where
\begin{align}\label{c3}
 C_3=\frac{i^{2h}2^{2h-1}\pi^{2h}}{(-1)^{v+1}k(2h)!} = -\frac{i^{2h}2^{2h-1}\pi^{2h}}{k(2h)!}
 \end{align}
 since $(-1)^{v+1}=-1$ as $r$ and $v$ are of same parity. So
\begin{align}\label{|L(2h,x|^2}
 \sum\limits_{\substack{\chi(\text{mod }  k) \\ \chi \text{ even}}}&|L(2h, \chi)|^2\\
 =&|C_3|^2\sum\limits_{\substack{\chi(\text{mod }  k) \\ \chi \text{ even}}}\left|\sum \limits_{q=0}^{2h-2}\binom{2h}{q}B_q \sum\limits_{j=1}^{k-1}\left(\frac{j}{k}\right)^{2h-q}\sum\limits_{m=1}^{k-1}\chi(m)e^{\frac{2\pi imj}{k}}\right|^2\nonumber\\
 =&|C_3|^2\sum\limits_{q_1=0}^{2h-2}\sum\limits_{q_2=0}^{2h-2}B_{q_1}B_{q_2}\binom{2h}{q_1}\binom{2h}{q_2}\sum\limits_{j=1}^{k-1}\left(\frac{j}{k}\right)^{2h-q_1}\sum\limits_{s=1}^{k-1}\left(\frac{s}{k}\right)^{2h-q_2}\nonumber\\ &\times\sum\limits_{m=1}^{k-1}\sum\limits_{n=1}^{k-1}e^{\frac{2\pi i(mj-ns)}{k}}\sum\limits_{\substack{\chi(\text{mod }  k) \\ \chi \text{ even}}}\chi(m)\overline{\chi(n)}.\nonumber
\end{align}
 Now  $\sum\limits_{\substack{\chi(\text{mod }  k) \\ \chi \text{ even}}}\chi(m)\overline{\chi(n)}=\sum\limits_{\substack{\chi(\text{mod }  k) \\ \chi \text{ even}}}\chi(mn^{-1})=\sum\limits_{\substack{\chi(\text{mod }  k) \\ \chi \text{ even}}}\chi(u)$ which is $\frac{\phi(k)}{2}$ if $u\equiv \pm1$(mod $k$) and $0$ otherwise for $k \geq 3$.
So we rewrite Equation (\ref{|L(2h,x|^2}) as
\begin{align}\label{first sum}
\sum\limits_{\substack{\chi(\text{mod }  k) \\ \chi \text{ even}}}&|L(2h, \chi)|^2\\
 =&|C_3|^2\frac{\phi(k)}{2}\sum\limits_{q_1=0}^{2h-2}\sum\limits_{q_2=0}^{2h-2}B_{q_1}B_{q_2}\frac{1}{k^{4h-q_1-q_2}}\binom{2h}{q_1}\binom{2h}{q_2}\nonumber\\
&\times\Bigg[\sum\limits_{\substack{1\leq m\leq k \\ (m, k)=1}}\left(\sum\limits_{j=1}^{k-1}j^{2h-q_1}e^{\frac{2\pi imj}{k}}\right)\left(\sum\limits_{s=1}^{k-1}s^{2h-q_2}e^{-\frac{2\pi ims}{k}}\right)\nonumber\\
&+\sum\limits_{\substack{1\leq m\leq k \\ (m, k)=1}}\left(\sum\limits_{j=1}^{k-1}j^{2h-q_1}e^{\frac{2\pi imj}{k}}\right)\left(\sum\limits_{s=1}^{k-1}s^{2h-q_2}e^{\frac{2\pi ims}{k}}\right)\Bigg].\nonumber
\end{align}
Now we use Equation (\ref{g(p)f(q)}) to replace\\ $\sum\limits_{\substack{1\leq m\leq k \\ (m, k)=1}}\left(\sum\limits_{j=1}^{k-1}j^{2h-q_1}e^{\frac{2\pi imj}{k}}\right)\left(\sum\limits_{s=1}^{k-1}s^{2h-q_2}e^{-\frac{2\pi ims}{k}}\right)$ with $-k^{2h-q_2}\sum\limits_{\substack{1\leq m\leq k \\ (m, k)=1}}\left(\sum\limits_{j=1}^{k-1}j^{2h-q_1}e^{\frac{2\pi imj}{k}}\right)+ \sum\limits_{a=0}^{2h-q_2-1}(-1)^{2h-q_2-a}\binom{2h-q_2}{a}k^{a}\sum\limits_{\substack{1\leq m\leq k \\ (m, k)=1}}\left(\sum\limits_{j=1}^{k-1}j^{2h-q_1}e^{\frac{2\pi imj}{k}}\right)\left(\sum\limits_{s=1}^{k-1}s^{2h-q_2-a}e^{\frac{2\pi ims}{k}}\right)$ to get the right-hand side of $\sum\limits_{\substack{\chi(\text{mod }  k) \\ \chi \text{ even}}}|L(2h, \chi)|^2$ above to be equal to
\begin{align*}
&\frac{\phi(k)}{2}\Bigg[
\sum\limits_{q_1=0}^{2h-2}\sum\limits_{q_2=0}^{2h-2}B_{q_1}B_{q_2}\frac{1}{k^{4h-q_1-q_2}}\binom{2h}{q_1}\binom{2h}{q_2}\times -k^{2h-q_2}\sum\limits_{\substack{1\leq m\leq k \\ (m, k)=1}}\left(\sum\limits_{j=1}^{k-1}j^{2h-q_1}e^{\frac{2\pi imj}{k}}\right)\\
&+\sum\limits_{q_1=0}^{2h-2}\sum\limits_{q_2=0}^{2h-2}B_{q_1}B_{q_2}\binom{2h}{q_1}\binom{2h}{q_2}\sum\limits_{a=0}^{2h-q_2-1}(-1)^{2h-q_2-a}
 \binom{2h-q_2}{a}k^{a-4h+q_1+q_2}\\&\times\sum\limits_{\substack{1\leq m\leq k \\ (m, k)=1}}\left(\sum\limits_{j=1}^{k-1}j^{2h-q_1}e^{\frac{2\pi imj}{k}}\right)\left(\sum\limits_{s=1}^{k-1}s^{2h-q_2-a}e^{\frac{2\pi ims}{k}}\right)
\\&+\sum\limits_{q_1=0}^{2h-2}\sum\limits_{q_2=0}^{2h-2}B_{q_1}B_{q_2}k^{-4h+q_1+q_2}\binom{2h}{q_1}\binom{2h}{q_2}\sum\limits_{\substack{1\leq m\leq k \\ (m, k)=1}}\left(\sum\limits_{j=1}^{k-1}j^{2h-q_1}e^{\frac{2\pi imj}{k}}\right)\left(\sum\limits_{s=1}^{k-1}s^{2h-q_2}e^{\frac{2\pi ims}{k}}\right)\Bigg].
\end{align*}

Now as in the case when  $r$ is odd, if we rearrange the summations in the sum
\begin{align*}
\Sigma'_0 :=\sum\limits_{q_1=0}^{2h-2}\sum\limits_{q_2=0}^{2h-2}B_{q_1}B_{q_2}\frac{1}{k^{4h-q_1-q_2}}\binom{2h}{q_1}\binom{2h}{q_2}\times -k^{2h-q_2}\sum\limits_{\substack{1\leq m\leq k \\ (m, k)=1}}\left(\sum\limits_{j=1}^{k-1}j^{2h-q_1}e^{\frac{2\pi imj}{k}}\right)
\end{align*}
in such a way that the innermost sum becomes $\sum\limits_{q_2=0}^{2h-2}B_{q_2}\binom{2h}{q_2}$, then by \cite[Theorem 12.15]{tom1976introduction}, since $2h\geq 4$, we get  $B_{2h}=\sum\limits_{q_2=0}^{2h}B_{q_2}\binom{2h}{q_2}=\sum\limits_{q_2=0}^{2h-2}B_{q_2}\binom{2h}{q_2}+B_{2h}$. Hence $\sum\limits_{q_2=0}^{2h-2}B_{q_2}\binom{2h}{q_2}=0$ and so $\Sigma'_0= 0$.\\
Hence Equation (\ref{first sum}) becomes
\begin{align*}
 \sum\limits_{\substack{\chi(\text{mod }  k) \\ \chi \text{ even}}}&|L(2h, \chi)|^2\\ = & |C_3|^2\frac{\phi(k)}{2}\Bigg[\sum\limits_{q_1=0}^{2h-2}\sum\limits_{q_2=0}^{2h-2}B_{q_1}B_{q_2}\binom{2h}{q_1}\binom{2h}{q_2}\sum\limits_{a=0}^{2h-q_2-1}(-1)^{2h-q_2-a}
 \binom{2h-q_2}{a}\nonumber\\&\times k^{a-4h+q_1+q_2}\sum\limits_{\substack{1\leq m\leq k \\ (m, k)=1}}\left(\sum\limits_{j=1}^{k-1}j^{2h-q_1}e^{\frac{2\pi imj}{k}}\right)\left(\sum\limits_{s=1}^{k-1}s^{2h-q_2-a}e^{\frac{2\pi ims}{k}}\right)
\nonumber\\&+\sum\limits_{q_1=0}^{2h-2}\sum\limits_{q_2=0}^{2h-2}B_{q_1}B_{q_2}k^{-4h+q_1+q_2}\binom{2h}{q_1}\binom{2h}{q_2}\sum\limits_{\substack{1\leq m\leq k \\ (m, k)=1}}\left(\sum\limits_{j=1}^{k-1}j^{2h-q_1}e^{\frac{2\pi imj}{k}}\right)\nonumber\\&\times\left(\sum\limits_{s=1}^{k-1}s^{2h-q_2}e^{\frac{2\pi ims}{k}}\right)\Bigg]\nonumber.
\end{align*}

In the above expression, since the left-hand side is a real number, we need to look only at the real part of the right-hand side. The imaginary part on the right-hand side is contributed by the sumations over $1\leq m\leq k, (m,k)=1$ in which the exponential function is present. For further computations, we consider the following  summations over $q_1$ and $q_2$ seperately in which the exponential function is present. Denote

\begin{align*}
 \Sigma_1':= &\sum\limits_{q_1=0}^{2h-2}\sum\limits_{q_2=0}^{2h-2}B_{q_1}B_{q_2}\binom{2h}{q_1}\binom{2h}{q_2}\sum\limits_{a=0}^{2h-q_2-1}(-1)^{2h-q_2-a}
 \binom{2h-q_2}{a}k^{a-4h+q_1+q_2}\\&\times\Re\left(\sum\limits_{\substack{1\leq m\leq k \\ (m, k)=1}}\left(\sum\limits_{j=1}^{k-1}j^{2h-q_1}e^{\frac{2\pi imj}{k}}\right)\left(\sum\limits_{s=1}^{k-1}s^{2h-q_2-a}e^{\frac{2\pi ims}{k}}\right)\right)\nonumber
\end{align*}
and
\begin{align*}
\Sigma_2' :=& \sum\limits_{q_1=0}^{2h-2}\sum\limits_{q_2=0}^{2h-2}B_{q_1}B_{q_2}k^{-4h+q_1+q_2}\binom{2h}{q_1}\binom{2h}{q_2}\\
 &\times\Re\left(\sum\limits_{\substack{1\leq m\leq k \\ (m, k)=1}}\left(\sum\limits_{j=1}^{k-1}j^{2h-q_1}e^{\frac{2\pi imj}{k}}\right)\left(\sum\limits_{s=1}^{k-1}s^{2h-q_2}e^{\frac{2\pi ims}{k}}\right)\right)\nonumber.
\end{align*}
After simplifying $\Sigma_1'$ and $\Sigma_2'$ using Lemma \ref{realjs} as in the case when $r$ is odd, we get
\begin{align*}
  \Sigma_1' =& \sum\limits_{q_1=0}^{2h-2}\sum\limits_{q_2=0}^{2h-2}\sum\limits_{a=0}^{2h-q_2-1}\sum\limits_{j=1}^{2h-q_1}\sum\limits_{s=1}^{2h-q_2-a}\sum\limits_{\alpha=1}^{2h+1-q_1-j}B_{q_1}B_{q_2}\binom{2h}{q_1}\binom{2h}{q_2}
 \binom{2h-q_2}{a}\\&\times
\binom{2h-q_1}{j}\binom{2h-q_2-a}{s}k^{a-4h+q_1+q_2+j+s} A_{2h-q_1-j, \alpha}\\&\times\Bigg[\Bigg(\sum\limits_{\substack{\beta=1\\ \alpha+\beta \text{ even}}}^{2h-q_2-a-s} \sum\limits_{c=0}^{\frac{\alpha+\beta}{2}}\sum\limits_{d=0}^{\frac{\alpha+\beta-2c}{2}}(-1)^{2h-q_2-a+c+d+\frac{\alpha+\beta}{2}}A_{2h-q_2-a-s, \beta}\nonumber\\
&\times\frac{(\alpha+\beta)(\alpha+\beta-c-1)!}{2^{2c+1}c!(\alpha+\beta-2c)!}
\binom{\frac{\alpha+\beta-2c}{2}}{d}SIN(\alpha+\beta-2d)\Bigg)\nonumber\\
& +\Bigg(\sum\limits_{\substack{\beta=1\\ \alpha+\beta \text{ odd}}}^{2h-q_2-a-s}\sum\limits_{e=0}^{\frac{\alpha+\beta-1}{2}}\sum\limits_{f=0}^{\frac{\alpha+\beta-2e-1}{2}}(-1)^{2h-q_2-a+e+f+\frac{\alpha+\beta+1}{2}}A_{2h-q_2-a-s, \beta}\nonumber\\
&\times\frac{(\alpha+\beta-e-1)!}{2^{2e+1}e!(\alpha+\beta-2e-1)!}\binom{\frac{\alpha+\beta-2e-1}{2}}{f}SIN(\alpha+\beta-2f-1)\Bigg)\Bigg]\nonumber
\end{align*}
and
\begin{align}\label{sigma2'}
       &\Sigma_2'= \sum\limits_{q_1=0}^{2h-2}\sum\limits_{q_2=0}^{2h-2}\sum\limits_{j=1}^{2h-q_1}\sum\limits_{s=1}^{2h-q_2}\sum\limits_{\alpha=1}^{2h-q_1-j+1}
  B_{q_1}B_{q_2}\binom{2h}{q_1}\binom{2h}{q_2}
\binom{2h-q_1}{j}\binom{2h-q_2}{s}\\
&\times k^{-4h+q_1+q_2+j+s} A_{2h+h_0-q_1-j, \alpha}\Bigg[\Bigg(\sum\limits_{\substack{\beta=1\\ \alpha+\beta \text{ even}}}^{2h-q_2-s+1} \sum\limits_{c=0}^{\frac{\alpha+\beta}{2}}\sum\limits_{d=0}^{\frac{\alpha+\beta-2c}{2}}(-1)^{c+d+\frac{\alpha+\beta}{2}}A_{2h-q_2-s, \beta}\nonumber\\
&\times\frac{(\alpha+\beta)(\alpha+\beta-c-1)!}{2^{2c+1}c!(\alpha+\beta-2c)!}
\binom{\frac{\alpha+\beta-2c}{2}}{d}SIN(\alpha+\beta-2d)\Bigg)\nonumber\\
& +\Bigg(\sum\limits_{\substack{\beta=1\\ \alpha+\beta \text{ odd}}}^{2h-q_2-s+1}\sum\limits_{e=0}^{\frac{\alpha+\beta-1}{2}}\sum\limits_{f=0}^{\frac{\alpha+\beta-2e-1}{2}}(-1)^{e+f+\frac{\alpha+\beta+1}{2}}A_{2h-q_2-s, \beta}\nonumber\\
&\times\frac{(\alpha+\beta-e-1)!}{2^{2e+1}e!(\alpha+\beta-2e-1)!}\binom{\frac{\alpha+\beta-2e-1}{2}}{f}SIN(\alpha+\beta-2f-1)\Bigg)\Bigg]\nonumber.
      \end{align}

In the above, $SIN(r)$ can be computed recusrively using Lemma \ref{lem:sineevenpowersum}. With the notations $\Sigma'_1$ and $\Sigma'_2$, we have

\begin{align}\label{L2h final}
 &\sum\limits_{\substack{\chi(\text{mod }  k) \\ \chi \text{ even}}}|L(2h, \chi)|^2= |C_3|^2\frac{\phi(k)}{2}\left(\Sigma'_1+\Sigma'_2\right).
\end{align}

Note that the sums $\Sigma_1'$ and $\Sigma_2'$  are almost same as the sums $\Sigma_1$ (Equation (\ref{sigma1})) and $ \Sigma_2$ (Equation (\ref{sigma2})) respectively. The difference between they are that
 where ever $2h+1$ is appearing in the  sum $\Sigma_i$, $2h$ is appearing in $\Sigma'_i$, and the upper limits of $q_1$ and $q_2$ in $\Sigma_i$ was $2h$ whereas in $\Sigma'_i$ the upper limit is $2h-2$.

 Recall that we obtained $\Sigma_1+\Sigma_2 = 0$ when $r$ is odd. Similar type of computations will give us here that $\Sigma_1'-\Sigma_2'=0$ and hence $\Sigma_1'=\Sigma_2'$.
Hence Equation (\ref{L2h final}) with the value of $C_3$ (Equation (\ref{c3})) substituted back gives us

\begin{theo}\label{identity:evencasefinalformula}
For $h \geq 2, h\in \N$, we have
\begin{align*}
 &\sum\limits_{\substack{\chi(\text{mod }  k) \\ \chi \text{ even}}}|L(2h, \chi)|^2 = \frac{2^{4h-2}\pi^{4h}}{((2h)!)^2k^2}\phi(k)\Sigma_2',\nonumber
\end{align*}
where $\Sigma_2'$ is as in Equation (\ref{sigma2'}).
\end{theo}

\subsection{Examples}
The mathematical software system SageMath \cite{sagemath} helped us to perform the computations in this section. We obtained the values of $SIN(2n)$ using the recursive procedure given in Lemma \ref{lem:sineevenpowersum} for $n=1,2,\ldots,6$ as
\begin{enumerate}
\item $\frac{1}{3} \, J_{2}(k)$, 
\item $\frac{1}{45} \, J_{4}(k)+\frac{2}{9} \, J_{2}(k)$, 
\item $\frac{2}{945} \, J_{6}(k)+ \frac{1}{45} \, J_{4}(k)+\frac{8}{45} \, J_{2}(k)$, 
\item $\frac{1}{4725} \, J_{8}(k)+ \frac{8}{2835} \, J_{6}(k)+ \frac{14}{675} \, J_{4}(k) + \frac{16}{105} \, J_{2} (k)$, 
\item $\frac{2}{93555} \, J_{10}(k)+ \frac{1}{2835} \, J_{8}(k)+ \frac{26}{8505} \, J_{6}(k)+ \frac{164}{8505} \, J_{4}(k) + \frac{128}{945} \, J_{2}(k)$,\\
and
\item $\frac{1382}{638512875} \, J_{12}(k)+\frac{4}{93555} \, J_{10}(k) + \frac{31}{70875} \, J_{8}(k) + \frac{556}{178605} \, J_{6}(k)  + \frac{3832}{212625} \, J_{4}(k) +\frac{256}{2079} \, J_{2}(k)$

\end{enumerate}
 respectively.

To evaluate  $\sum\limits_{\substack{\chi(\text{mod }  k) \\ \chi \text{ odd}}}|L(5, \chi)|^2 $, we evaluated $\Sigma_2$ appearing in Theorem \ref{identity:oddcasefinalformula}  and obtained $\Sigma_2=-\frac{5}{16632 \, k^{8}} (J_{10}(k)- 22 \, J_{4}(k) - 231 \, J_{2}(k) )$. Thus, by Theorem \ref{identity:oddcasefinalformula},
\begin{align*}
\sum\limits_{\substack{\chi(\text{mod }  k) \\ \chi \text{ odd}}}|L(5, \chi)|^2=\frac{\pi^{10}\phi(k)}{187110 \, k^{10}}\left( J_{10}(k) - 22 \,  J_{4}(k) - 231 \, J_{2}(k)\right).
\end{align*}
This agrees with the value of $\sum\limits_{\substack{\chi(\text{mod }  k) \\ \chi \text{ odd}}}|L(5, \chi)|^2$  computed by X. Lin in \cite{lin2019mean}.

To compute $\sum\limits_{\substack{\chi(\text{mod }  k) \\ \chi \text{ even}}}|L(6, \chi)|^2 $, we evaluated   $\Sigma_2'$ appearing in Theorem \ref{identity:evencasefinalformula} and obtained   $\Sigma_2'=\frac{1 }{2522520 \, k^{10}} ( 691 \, J_{12} (k) + 2860 \, J_{6} (k)+ 63063 \, J_{4} (k)+ 573300 \, J_{2} (k)) $. Thus by Theorem \ref{identity:evencasefinalformula}
\begin{align*}
\sum\limits_{\substack{\chi(\text{mod }  k) \\ \chi \text{ even}}}|L(6, \chi)|^2  =&\frac{ \pi^{12} \phi(k)}{1277025750 \, k^{12}} (691 \,  J_{12}(k)+ 2860 \,  J_{6}(k)\\&+ 63063 \, J_{4}(k) + 573300 \,  J_{2} (k)).
\end{align*}

Note that the values of $\sum\limits_{\substack{\chi(\text{mod }  k) \\ \chi(-1)=(-1)^r}}|L(r, \chi)|^2$ for $r=1,2$ appear in \cite{alkan2011mean} and for $r=3,4$  in \cite{thomas2023mean}.

\section{Acknowledgements}
 The first author thanks the Kerala State Council for Science,Technology and Environment, Thiruvananthapuram, Kerala, India for providing financial support for carrying out this research work.

\bibliographystyle{plain}

\end{document}